\documentclass[12]{amsart}
\usepackage{latexsym}
\usepackage{amsthm}
\usepackage{amssymb}
\usepackage{amsfonts, amsmath}

\newtheorem{Prop}{Proposition}
\newtheorem*{Propo}{Proposition}

\newtheorem{lm}[Prop]{Lemma}
\newtheorem*{Thm}{Theorem}

\title {\title\\ On the Integral homology and counterexamples\\ to the
Andreotti-Grauert conjecture}
\author{\bf MOHAMMED MOU\c{C}OUF AND YOUSSEF ALAOUI} \subjclass{32E10, 32E40.}
\keywords{$q$-complete spaces; $q$-Runge domains; $q$-convex
functions; cohomologically $q$-complete spaces.}
\date{}
\begin{document}
 \maketitle
 \begin{abstract}In this paper, we prove by means of a counterexample that there exist pair of integers
$(n,p)$ with $n\geq 3$, $1<p\leq n-1$, and open sets $D$ in $\mathbb{C}^{n}$ which are cohomologically $p$-complete with respect to the structure sheaf ${\mathcal{O}}_{D}$
such that the cohomology group $H_{n+p}(D, \mathbb{Z})$ does not vanish. In particular $D$ is not $p$-complete.
\end{abstract}
\maketitle \setcounter{page}{1}\noindent
\section{Introduction} By
the theory of Andreotti and Grauert ~\cite{ref2} it is known that
a $q$-complete complex space is always cohomologically
$q$-complete.\\
\hspace*{.1in}It was shown in ~\cite{ref7} that if $X$ is a Stein manifold and $\Omega$ is a cohomologically $q$-complete open set in $X$ with respect to the structure sheaf ${\mathcal{O}}_{\Omega}$, then $\Omega$ is $q$-complete, if it has a smooth boundary.\\
\hspace*{.1in}If now $X$ is a $q$-complete space of complex dimension
$n$, then $H_{p}(X,\mathbb{Z})=0$ for $p\geq n+q$ and
$H_{n+q-1}(X,\mathbb{Z})$ is free (see ~\cite{ref6} and
~\cite{ref12}). Considering the short exact sequence given by the
universal coefficient theorem
$$0\rightarrow
Ext_{\mathbb{Z}}(H_{p-1}(X,\mathbb{Z}),\mathbb{Z})\rightarrow
H^{p}(X,\mathbb{Z})\rightarrow
Hom_{\mathbb{Z}}(H_{p}(X,\mathbb{Z}),\mathbb{Z})\rightarrow 0$$
we see that $H^{p}(X,\mathbb{Z})=0$ for $p\geq n+q$.\\
\hspace*{.1in}If $X$ is a cohomologically $q$-complete complex
space of dimension $n$, it follows from ~\cite{ref10} that
$H_{p}(X,\mathbb{C})=0$ for all integers $p\geq n+q$. Since
$H^{p}(X,\mathbb{C})$, for an $n$-dimensional complex manifold $X$,
are Frechet spaces whose topological duals are isomorphic to
$H_{p}(X,\mathbb{C})$ (see ~\cite{ref11}), then $X$ satisfies the
topological property $H^{p}(X,\mathbb{C})=0$ for $p\geq n+q$.\\
\hspace*{.1in}However, it seems unknown if $H^{p}(X,\mathbb{Z})$,
$p\geq n+q$, vanishes, if $X$ is assumed to be only cohomologically $q$-complete.\\
\hspace*{.1in}In this article, we prove that there exist pair of integers $(n,p)$
with $1<p\leq n-1$, and cohomologically $p$-complete open sets $D$
with respect to ${\mathcal{O}}_{D}$
in $\mathbb{C}^{n}$ such that $H_{n+p}(D,\mathbb{Z})\neq
0$. It is clear that $D$ is not $p$-complete.
\section{Preliminaries}
Let $\Omega$ be a an open set in $\mathbb{C}^{n}$ with complex
coordinates $z_{1},\cdots,z_{n}$. Then it is known that a function
$\phi\in C^{\infty}(\Omega)$ is said to be $q$-convex if the
hermetian form
$L_{z}(\phi,\xi)=\displaystyle\sum_{i,j}\frac{\partial^{2}\phi(z)}{\partial{z}_{i}\partial{\overline{z}}_{j}}\xi_{i}\overline{\xi}_{j}$
has at least $n-q+1$ positive eingenvalues at each point $z\in \Omega.$\\
\hspace*{.1in}A function $\rho\in C^{o}(\Omega,\mathbb{R})$ is
called $q$-convex with corners if $\rho$ is locally a maximum of
finite number
of $q$-convex functions.\\
\hspace*{.1in}We say that $\Omega$ is $q$-complete (resp.
$q$-complete with corners) if there exists an exhaustion function
$\phi$ on $\Omega$
which is $q$-convex (resp. $q$-convex with corners).\\
\hspace*{.1in}An open subset $D$ of $\Omega$ is called $q$-Runge if
for every compact set $K\subset D$, there is a $q$-convex exhaustion
function $\phi\in C^{\infty}(\Omega)$ such that
$$K\subset\{x\in\Omega: \phi(x)<0\}\subset\subset D$$
It is shown in ~\cite{ref2} that if $\Omega$ is $q$-complete, then
$\Omega$ is cohomologically $q$-complete, which means that for every
coherent analytic sheaf ${\mathcal{F}}$ on $\Omega$ and $p\geq q$,
the cohomology groups $H^{p}(\Omega,{\mathcal{F}})$ vanish.
Moreover, if $D$ is $q$-Runge in $\Omega$, then for every
${\mathcal{F}}\in coh(\Omega)$ the restriction map
$$H^{p}(\Omega,{\mathcal{F}})\rightarrow H^{p}(D,{\mathcal{F}})$$
has a dense image for all $p\geq q-1.$\\
\\
\hspace*{.1in}The purpose of the present article is to prove the
following theorem
\begin{Thm}Let $(n,q)$ be a pair of integers such that $n\geq 3,$
$1<q<n$, and $q$ does not divide $n$.
We put $m=[\frac{n}{q}]$ and $r=n-mq$.
Then there exists a
cohomologically $(\tilde{q}-3)$-complete domain $D$ in
$\mathbb{C}^{n}$ such that $H_{n+\tilde{q}-3}(D,\mathbb{Z})\neq
0$, if $r=1$ and $m>q\geq 2$.
Here $\tilde{q}=n-[\frac{n}{q}]+1$ and $[\frac{n}{q}]$ is the
integral part of $\frac{n}{q}$.
\end{Thm}
\section{Proof of the theorem}
Let $(n,q)$ be a pair of integers with $n\geq 3$ and
$1\leq q\leq n.$ We put $m=[\frac{n}{q}]$ and suppose that
$r=n-mq>0.$ We consider the functions $\phi_{1},\cdots, \phi_{m+1}$
defined on $\mathbb{C}^{n}$ by
$$\phi_{j}=\sigma_{j}+\displaystyle\sum_{i=1}^{m}\sigma_{i}^{2}
-\frac{1}{4}||z||^{2}+N||z||^{4}, \ \ j=1,\cdots,m,$$ and
$$\phi_{m+1}=-\sigma_{1}-\cdots-\sigma_{m}+\displaystyle\sum_{i=1}^{m}\sigma_{i}^{2}
-\frac{1}{4}||z||^{2}+N||z||^{4},$$ where
$\sigma_{j}=y_{j}+\displaystyle\sum_{i=m+1}^{n}|z_{i}|^{2}-
(m+1)(\displaystyle\sum_{i=m+(j-1)(q-1)+1}^{m+j(q-1)}|z_{i}|^{2})$
for $j=1,\cdots,m$, $z_{j}=x_{j}+iy_{j}$. Then it is known from
~\cite{ref5} that all $\phi_{j}$, $1\leq j\leq m+1$, are $q$-convex
on $\mathbb{C}^{n}$, if $N>0$ is sufficiently large and, if
$\rho=Max\{\phi_{j}, 1\leq j\leq m+1\}$, then there exists a small
constant $\varepsilon_{o}>0$ such that the set
$D_{\varepsilon_{o}}=\{z\in\mathbb{C}^{n}:
\rho(z)<-\varepsilon_{o}\}$ is relatively compact in the unit ball
$B=B(0,1)$, if $N$ is large enough.\\
\hspace*{.1in}Let now $\varepsilon>\varepsilon_{0}$ and choose Stein
open sets $U_{1},\cdots,U_{k}\subset\subset D_{\varepsilon_{o}}$
covering $\partial{D_{\varepsilon}}$ and functions $\theta_{j}\in
C^{\infty}_{o}(U_{j}, \mathbb{R}^{+})$ such that
$\displaystyle\sum_{j=1}^{k}\theta_{j}(x)>0$ at any point $x\in
\partial{D_{\varepsilon}}$. There exist sufficiently small constants
$c_{1}>0,\cdots, c_{k}>0$ such that the functions
$\phi_{i,j}=\phi_{i}-\displaystyle\sum_{i=1}^{j}c_{i}\theta_{i}$
are $q$-convex for $i=1,\cdots, m+1$ and $1\leq j\leq k$. We
define $\phi_{i,0}=\phi_{i}$ for $i=1,\cdots,m+1,$
$D_{o}=D_{\varepsilon}$ and $D_{j}=\{z\in D_{\varepsilon_{o}}:
\rho_{j}(z)<-\varepsilon\}$, where
$$\rho_{j}=\rho-\displaystyle\sum_{i=1}^{j}c_{i}\theta_{i}, \ \ j=1,\cdots, k$$
Then $\rho_{j}$ are $q$-convex with corners, $D_{o}\subset
D_{1}\subset\cdots\subset D_{k}$,\\ $D_{o}\subset\subset
D_{k}\subset\subset D_{\varepsilon_{o}}$ and $D_{j}\setminus
D_{j-1}\subset\subset U_{j}$ for $j=1,\cdots, k$.
\begin{lm}Let ${\mathcal{F}}$ be a coherent analytic sheaf on $D_{\varepsilon_{0}}$.
Then the restriction map  $H^{p}(D_{j+1}, {\mathcal{F}})\rightarrow H^{p}(D_{j}, {\mathcal{F}})$
is surjective for every $p\geq \tilde{q}-2$ and all $0\leq j\leq
k-1.$ In particular, $dim_{\mathbb{C}}H^{p}(D_{j},{\mathcal{F}})<\infty$
for $p\geq \tilde{q}-2$ and $ j=0,\cdots, k$.
\end{lm}
\begin{proof}
\hspace*{.1in}We first prove that for every $p\geq \tilde{q}-2$,
$H^{p}(D_{j}\cap U_{l},{\mathcal{F}})=0$ for all $0\leq j\leq k$
and $1\leq l\leq k$. We fix $j\in\{0,\cdots,k\}$ and we write
$D_{j}\cap U_{l}=D'_{1}\cap\cdots\cap D'_{m+1},$ where
$D'_{i}=\{z\in U_{l}: \phi_{i,j}(z)<-\varepsilon\}$ are clearly
$q$-complete and $q$-Runge in $U_{l}$. Then for any integer $t\leq
m$, $D'_{i_{1}}\cap \cdots\cap D'_{i_{t}}$ is
$(\tilde{q}-1)$-Runge in $U_{l}$ for all
$i_{1},\cdots,i_{t}\in\{1,\cdots,m+1\}$. In fact, let $K\subset
D'_{i_{1}}\cap\cdots\cap D'_{i_{t}}$ be an arbitrary compact
subset. There exists for any $i\in\{i_{1},\cdots,i_{t}\}$ a
$q$-convex exhaustion function $\psi_{i}$ on $U_{l}$ such that
$$K\subset \{x\in U_{l}: \psi_{i}(x)<0\}\subset\subset D'_{i}$$
Define $\psi=Max(\psi_{i_{1}},\cdots,\psi_{i_{t}})$. Then
$K\subset \{x\in U_{l}: \psi(x)<0\}\subset\subset
D'_{i_{1}}\cap\cdots\cap D'_{i_{t}}$ and, $\psi$ can be
approximated in the $C^{O}$-topology by smooth $(t(q-1)+1)$-convex
functions (see ~\cite{ref11}). Since $q$ does not divide $n$ and
$t\leq m$, then $t(q-1)+1\leq \tilde{q}-1$. Therefore a suitable
smooth $(\tilde{q}-1)$-convex approximation of $\psi$ shows that
$D'_{i_{1}}\cap \cdots\cap D'_{i_{t}}$ is $(\tilde{q}-1)$-Runge in
$U_{l}$. This implies that for $r\geq \tilde{q}-2$, the image of
$H^{r}(U_{l},{\mathcal{F}})$ is dense in the cohomology group
$H^{r}(D'_{i_{1}}\cap\cdots\cap D'_{i_{t}},{\mathcal{F}})$ which
is separated, since $\tilde{q}-2\geq 2$ ($n\geq 4$ and $q\nmid
n$). Thus $H^{r}(D'_{i_{1}}\cap\cdots\cap
D'_{i_{t}},{\mathcal{F}})=0$ for all $r\geq \tilde{q}-2$. Then, by
(~\cite{ref9}, Proposition $1$), we have
$$H^{r}(D_{j}\cap U_{l},{\mathcal{F}})\cong H^{r+m}(D'_{1}\cup\cdots\cup
D'_{m+1},{\mathcal{F}}) \ \ if \ \ r\geq \tilde{q}-2$$ Since
$U_{l}\backslash D'_{1}\cup\cdots\cup D'_{m+1}$ has no compact
connected components, then $D'_{1}\cup\cdots\cup D'_{m+1}$ is
$n$-Runge in $U_{l}$ (see ~\cite{ref4}). But for $r\geq
\tilde{q}-2$, $r+m\geq n-1$ and hence
$H^{r+m}(D'_{1}\cup\cdots\cup D'_{m+1},{\mathcal{F}})=0.$\\
\hspace*{.1in}Now since $D_{j+1}=D_{j}\cup (D_{j+1}\cap U_{j+1}),$
then the Mayer-Vietoris sequence for cohomology
$$\cdots\rightarrow H^{p}(D_{j+1},{\mathcal{F}})\rightarrow
H^{p}(D_{j},{\mathcal{F}})\oplus H^{p}(D_{j+1}\cap
U_{j+1},{\mathcal{F}})\rightarrow H^{p}(D_{j}\cap
U_{j+1},{\mathcal{F}})\rightarrow\cdots$$ implies that the
restriction map $H^{p}(D_{j+1},{\mathcal{F}})\rightarrow
H^{p}(D_{j},{\mathcal{F}})$ is surjective for all $p\geq
\tilde{q}-2.$
\end{proof}
We now choose $n$ and $q$ such that $m=[\frac{n}{q}]\geq q$ and $r=n-mq\geq 1$, and
let ${\mathcal{O}}$ be the sheaf of germs of holomorphic functions on $B$.
Then we have the following
\begin{lm}The restriction map $H^{p}(D_{\varepsilon_{o}}, {\mathcal{O}})\rightarrow H^{p}(D_{\varepsilon}, {\mathcal{O}})$ has dense image for every $p\geq \tilde{q}-3$ and $\varepsilon\geq
\varepsilon_{o}$
\end{lm}
\begin{proof}Let $T$ be the set of all real numbers
$\varepsilon\geq \varepsilon_{o}$ such that
$H^{p}(D_{\varepsilon}, {\mathcal{O}})\rightarrow
H^{p}(D_{\varepsilon_{1}}, {\mathcal{O}})$ has a dense image for
every real number $\varepsilon_{1}>\varepsilon$ and all $p\geq
\tilde{q}-3.$ Then $T\neq\emptyset$. In fact, choose
$\varepsilon>\varepsilon_{o}$ such that
$-\varepsilon<Min_{\overline{B}\setminus
D_{\varepsilon_{o}}}\{\phi_{i}(z), i=1,\cdots, m+1\}$, and let
$\varepsilon_{1}>\varepsilon$. If
$D_{\varepsilon_{1}}\neq\emptyset,$ then $D_{i}=\{z\in B:
\phi_{i}(z)<-\varepsilon_{1}\}$ and $D'_{i}=\{z\in B:
\phi_{i}(z)<-\varepsilon\}$ are relatively compact in
$D_{\varepsilon_{o}}$, $q$-complete and $q$-Runge in $B$.\\Note
first that for every $i_{1},\cdots,i_{m-1}\in \{1,\cdots,m+1\}$,
the cohomology group $H^{p}(D_{i_{1}}\cap \cdots \cap
D_{i_{m-1}},\mathcal{O})=0$ for $p\geq \tilde{q}-2$, since
$D_{i_{1}}\cap \cdots \cap D_{i_{m-1}}$ is
$((m-1)(q-1)+1)$-complete and $\tilde{q}-2\geq ((m-1)(q-1)+1)$.
Next, we show, exactly as in the proof of lemma $1$, that for
every $i_{1},\cdots,i_{m}\in \{1,\cdots,m+1\}$, $D_{i_{1}}\cap
\cdots \cap D_{i_{m}}$ is $(\tilde{q}-1)$-Runge in $B$, which
implies that the restriction map $H^{p}(D_{i_{1}}\cap \cdots \cap
D_{i_{m-1}},\mathcal{O})\rightarrow H^{p}(D_{i_{1}}\cap \cdots
\cap D_{i_{m}},\mathcal{O})$ has a dense image for all $p\geq
\tilde{q}-2$. This shows that $H^{p}(D_{i_{1}}\cap \cdots \cap
D_{i_{m}},\mathcal{O})=0$ for $p\geq \tilde{q}-2$. Then, by using
(~\cite{ref9}, Proposition $1$), we obtain
$$H^{p}(D_{\varepsilon_{1}},\mathcal{O})\cong H^{p+m}(D_{1}\cup
\cdots \cup D_{m+1},\mathcal{O}) \mbox{\;for\;} p\geq\tilde{q}-2.$$
Similarly
$$H^{p}(D_{\varepsilon},\mathcal{O})\cong H^{p+m}(D'_{1}\cup
\cdots \cup D'_{m+1},\mathcal{O}) \mbox{\;for\;}
p\geq\tilde{q}-2.$$ Since $D_{1}\cup \cdots \cup D_{m+1}$ is
$n$-Runge in $B$ and contained in the open set $D'_{1}\cup \cdots
\cup D'_{m+1}\subset\subset D_{\varepsilon_{0}}$, then the
restriction map $$H^{p}(D'_{1}\cup \cdots \cup
D'_{m+1},\mathcal{O})\rightarrow H^{p}(D_{1}\cup \cdots \cup
D_{m+1},\mathcal{O})$$ has a dense image for $p\geq n-1$, which
means that $H^{p}(D_{\varepsilon},\mathcal{O})\rightarrow
H^{p}(D_{\varepsilon_{1}},\mathcal{O})$ has a dense image for
$p\geq \tilde{q}-2$.  We are now going to show that
$H^{\tilde{q}-3}(D_{\varepsilon},{\mathcal{O}})=H^{\tilde{q}-3}(D_{\varepsilon_{1}},{\mathcal{O}})=0$.
To see this, let $\Omega_{m}=\{z\in D_{1}\cap\cdots\cap D_{m}: \phi_{m+1}(z)>-\varepsilon\}$
and $S_{m}=\{z\in D_{1}\cap\cdots\cap D_{m}: \phi_{m+1}(z)\leq -\varepsilon\}$.
Then $H_{S_{m}}^{p}(D_{1}\cap\cdots\cap D_{m},{\mathcal{O}})=0$ for all $p\leq n-q$, where
$H_{S_{m}}^{j}(D_{1}\cap\cdots\cap D_{m},{\mathcal{O}})$ is the j-th group of cohomology of
$D_{1}\cap\cdots\cap D_{m}$ with support in $S_{m}$. In fact, for each point $\xi\in D_{1}\cap\cdots\cap D_{m}$
there exists ~\cite{ref2} a fundamental system of connected Stein neighborhoods $U\subset D_{1}\cap\cdots\cap D_{m}$ of $\xi$ such that $H^{j}(U\cap\Omega_{m},{\mathcal{O}})=0$ for $0<j<n-q$ and,
the restriction map
$$\Gamma(U,{\mathcal{O}})\longrightarrow \Gamma(U\cap \Omega_{m},{\mathcal{O}})$$
is an isomorphism. It follows from ~\cite{ref8} that $\underline{H^{j}_{S_{m}}}(\mathcal{O})=0$ for $0\leq j\leq n-q,$
where $\underline{H^{j}_{S_{m}}}(\mathcal{O})$ is the cohomology sheaf of $D_{1}\cap\cdots\cap D_{m}$ with coefficient in ${\mathcal{O}}$ and support in $S_{m}$. By (~\cite{ref8}) there is a spectral sequence
$$H^{p}_{S_{m}}(D_{1}\cap\cdots\cap D_{m},{\mathcal{O}})\Longleftarrow E_{2}^{p,q}=H^{p}(D_{1}\cap\cdots\cap D_{m},\underline{H^{q}_{S_{m}}}({\mathcal{O}})) $$
Since $\underline{H^{p}_{S_{m}}}({\mathcal{O}})=0$ for $p\leq n-q$, then the abutment
$H^{p}_{S_{m}}(D_{1}\cap\cdots\cap D_{m},{\mathcal{O}})=0$ for $p\leq n-q$.\\
\hspace*{.1in}Now it follows from the exact sequence of local cohomology
$$\cdots\rightarrow H^{p}_{S_{m}}(D_{1}\cap\cdots\cap D_{m},{\mathcal{O}})\rightarrow H^{p}(D_{1}\cap\cdots\cap D_{m},{\mathcal{O}})\rightarrow H^{p}(\Omega_{m},{\mathcal{O}})\rightarrow
H^{p+1}_{S_{m}}(D_{1}\cap\cdots\cap D_{m},{\mathcal{O}})\rightarrow\cdots$$
that $H^{p}(D_{1}\cap\cdots\cap D_{m},{\mathcal{O}})\cong H^{p}(\Omega_{m},{\mathcal{O}})$
for $0\leq p\leq n-q-1$.\\
\hspace*{.1in}Note that since $\Omega_{m}=\{z\in (D_{1}\cap\cdots\cap D_{m})\cup D_{m+1}: \phi_{m+1}(z)>-\varepsilon\}$, one can verify exactly as for $H_{S_{m}}^{p}(D_{1}\cap\cdots\cap D_{m},{\mathcal{O}})$ that, if\\ $S'_{m}=\{z\in (D_{1}\cap\cdots\cap D_{m})\cup D_{m+1}: \phi_{m+1}(z)\leq -\varepsilon\}=D_{m+1}\cup S_{m}$, then $H_{S'_{m}}^{p}((D_{1}\cap\cdots\cap D_{m})\cup D_{m+1},{\mathcal{O}})=0$ for $p\leq n-q$ and, therefore
$$H^{p}((D_{1}\cap\cdots\cap D_{m})\cup D_{m+1},{\mathcal{O}})\cong H^{p}((D_{1}\cap\cdots\cap D_{m})\cup D_{m+1}\setminus S'_{m},{\mathcal{O}})=H^{p}(\Omega_{m},{\mathcal{O}})\cong H^{p}(\Omega_{m}\cup D_{m+1},{\mathcal{O}})$$
for $q\leq p\leq n-q-1$.\\
\hspace*{.1in}Now, since, in addition,
$\Omega_{m}\cap D_{m+1}=\emptyset$ , then, by using the Mayer-Vietoris sequence for cohomology
\begin{center}
$\cdots\rightarrow H^{p}((D_{1}\cap\cdots\cap D_{m})\cup D_{m+1},{\mathcal{O}})\rightarrow
H^{p}((D_{1}\cap\cdots\cap D_{m}),{\mathcal{O}})\oplus H^{p}(D_{m+1},{\mathcal{O}})
\rightarrow H^{p}(D_{\varepsilon},{\mathcal{O}})\rightarrow H^{p+1}((D_{1}\cap\cdots\cap D_{m})\cup D_{m+1},{\mathcal{O}})\rightarrow H^{p+1}((D_{1}\cap\cdots\cap D_{m}),{\mathcal{O}})\oplus H^{p+1}(D_{m+1},{\mathcal{O}})\rightarrow\cdots$
\end{center}
and the fact that $D_{m+1}$ is $q$-complete, we find
that $H^{p}(D_{\varepsilon},{\mathcal{O}})=0$ for $q\leq p\leq n-q-2$. This implies that
$H^{\tilde{q}-3}(D_{\varepsilon},{\mathcal{O}})=0$ and similarly
$H^{\tilde{q}-3}(D_{\varepsilon_{1}},{\mathcal{O}})=0$,
which proves that $\varepsilon\in T$.\\
\hspace*{.1in}To see that $T$ is open in
$[\varepsilon_{o},+\infty[$, it is sufficient to prove that if
$\varepsilon\in T$, $\varepsilon>\varepsilon_{o},$ then there is
$\varepsilon_{o}<\varepsilon'<\varepsilon$ such that
$\varepsilon'\in T$. For this, we consider as in lemma $1$, finitely
many Stein open sets $U_{i}\subset\subset D_{\varepsilon_{o}},$
$i=1,\cdots, k,$ such that
$\partial{D}_{\varepsilon}\subset\displaystyle\bigcup_{i=1}^{k}U_{i}$
and functions $\theta_{j}\in C^{\infty}(U_{j},\mathbb{R}^{+})$ with
compact support such that
$\displaystyle\sum_{j=1}^{k}\theta_{j}(x)>0$ at any point
$x\in\partial{D}_{\varepsilon}$. Next we define
$D_{j}(\varepsilon)=\{z\in D_{\varepsilon_{o}}:
\rho_{j}(z)<-\varepsilon\}$, where
$\rho_{j}(z)=Max(\phi_{1}-\displaystyle\sum_{i=1}^{j}c_{i}\theta_{i},\cdots,\phi_{m+1}-\displaystyle\sum_{i=1}^{j}c_{i}\theta_{i})$
with $c_{i}>0$ sufficiently small so that the functions
$\phi_{i}-\displaystyle\sum_{i=1}^{j}c_{i}\theta_{i}$ are $q$-convex
for $1\leq i\leq m+1$ and $1\leq j\leq k$. Then, by lemma $1$, the
restriction map $H^{p}(D_{k}(\varepsilon),{\mathcal{O}})\rightarrow
H^{p}(D_{\varepsilon},{\mathcal{O}})$ is surjective for $p\geq
\tilde{q}-2$ and, there exists
$\varepsilon_{o}<\varepsilon'<\varepsilon$ such that
$D_{\varepsilon}\subset D_{\varepsilon'}\subset D_{k}(\varepsilon)$.
If now $\varepsilon'\leq \alpha\leq \varepsilon$, then we have
$$D_{\alpha}\subset D_{\varepsilon'}\subset D_{k}(\varepsilon)\subset
D_{k}(\alpha)\subset\subset D_{\varepsilon_{o}}.$$ Since, by lemma
$1$, the restriction map
$H^{p}(D_{k}(\alpha),{\mathcal{O}})\rightarrow
H^{p}(D_{\alpha},{\mathcal{O}})$ is surjective for $p\geq
\tilde{q}-2$ and $H^{\tilde{q}-3}(D_{\alpha},{\mathcal{O}})=0$, then
$H^{p}(D_{\varepsilon'},{\mathcal{O}})\rightarrow
H^{p}(D_{\alpha},{\mathcal{O}})$ is surjective for $p\geq
\tilde{q}-3$. For $\alpha>\varepsilon$, we have $D_{\alpha}\subset
D_{\varepsilon}\subset D_{\varepsilon'}\subset
D_{k}(\varepsilon).$ Since
$H^{p}(D_{k}(\varepsilon),{\mathcal{O}})\rightarrow
H^{p}(D_{\varepsilon},{\mathcal{O}})$ is surjective for $p\geq
\tilde{q}-3$, and $H^{p}(D_{\varepsilon},{\mathcal{O}})\rightarrow
H^{p}(D_{\alpha},{\mathcal{O}})$ has a dense image for $p\geq
\tilde{q}-3$, then $H^{p}(D_{\varepsilon'},\mathcal{O})\rightarrow
H^{p}(D_{\alpha},\mathcal{O})$ has a dense image for $p\geq \tilde{q}-3$,
which implies that $\varepsilon'\in T$.\\
\hspace*{.1in}In order to prove that $T$ is closed, we consider a
sequence of real numbers $\varepsilon_{j}\in T$, $j\geq 0$, such
that $\varepsilon_{j}\searrow\varepsilon$ and a Stein open
covering ${\mathcal{U}}=(U_{i})_{i\in I}$ of $D_{\varepsilon_{o}}$
with a countable base of open subsets of $D_{\varepsilon_{o}}$. We
fix $p\geq \tilde{q}-3$. Then the restriction map of spaces of
cocycles
$Z^{p}({\mathcal{U}}|_{D_{\varepsilon_{j+1}}},{\mathcal{O}})
\rightarrow
Z^{p}({\mathcal{U}}|_{D_{\varepsilon_{j}}},{\mathcal{O}})$ has a
dense image for $j\geq 0$. Therefore, by (~\cite{ref2}, p. $246$),
the restriction map
$Z^{p}({\mathcal{U}}|_{D_{\varepsilon}},{\mathcal{O}}) \rightarrow
Z^{p}({\mathcal{U}}|_{D_{\varepsilon_{j}}},{\mathcal{O}})$ has
also a dense image. Let $\varepsilon'>\varepsilon$ and $j\in
\mathbb{N}$ such that $\varepsilon'>\varepsilon_{j}$.  Since
$\varepsilon_{j}\in T$, then
$Z^{p}({\mathcal{U}}|_{D_{\varepsilon_{j}}},{\mathcal{O}})
\rightarrow
Z^{p}({\mathcal{U}}|_{D_{\varepsilon'}},{\mathcal{O}})$ has a
dense image, and hence $\varepsilon\in T.$ This completes the
proof of lemma $3$.
\end{proof}
Let now $A$ be the set of all real numbers $\varepsilon\geq \varepsilon_{0}$ such that
$H^{p}(D_{\varepsilon},{\mathcal{O}})=0$ for all $p\geq \tilde{q}-2$.
Then in the situation described above we have
\begin{lm}The set $A$ is not empty and, for every $\varepsilon\in A$ with $\varepsilon>\varepsilon_{0}$,
there is $\varepsilon_{o}\leq \varepsilon'<\varepsilon$ such that
$\varepsilon'\in A$.
\end{lm}
\begin{proof} Let $\varepsilon_{1}>\varepsilon_{o}$ be such that
$$-\varepsilon_{1}< Inf_{z\in \partial{D_{\varepsilon_{o}}}}\{\phi_{i}(z), i=1, 2,\cdots, m+1\}.$$
Then $[\varepsilon_{1},+\infty[\subset A.$ In fact, let
$\varepsilon\geq \varepsilon_{1}$ and write
$D_{\varepsilon}=D_{\varepsilon,1}\cap
D_{\varepsilon,2}\cap\cdots\cap D_{\varepsilon,m+1},$ where
$D_{\varepsilon,i}=\{z\in B:
\phi_{i}(z)<-\varepsilon\}\subset\subset D_{\varepsilon_{0}}$ is $q$-complete and $q$-Runge in $B$ for all $i\in\{1,\cdots, m+1\}$.
Moreover, for every integers $t\leq m$ and $i_{1},\cdots, i_{t}\in\{1,\cdots,m+1\}$,
$D_{\varepsilon,i_{1}}\cap
D_{\varepsilon,i_{2}}\cap\cdots\cap D_{\varepsilon,i_{t}}$ is $(tq-t+1)$-Runge in $B$ and $\tilde{q}-2\geq tq-t$, then $D_{\varepsilon,i_{1}}\cap
D_{\varepsilon,i_{2}}\cap\cdots\cap D_{\varepsilon,i_{t}}$ is cohomologically $(\tilde{q}-2)$-complete.
Therefore, by (~\cite{ref9}, Proposition $1$)
$$H^{r}(D_{\varepsilon},{\mathcal{O}})\cong H^{r+m}(D_{\varepsilon,1}\cup
D_{\varepsilon,2}\cup\cdots\cup D_{\varepsilon,m+1},{\mathcal{O}}), \ \ for \ \ r\geq \tilde{q}-2$$
Since $D_{\varepsilon,1}\cup
D_{\varepsilon,2}\cup\cdots\cup D_{\varepsilon,m+1}$ is $n$-Runge in $B$, then\\
$H^{p}(B,{\mathcal{O}})\rightarrow H^{p}(D_{\varepsilon,1}\cup
D_{\varepsilon,2}\cup\cdots\cup D_{\varepsilon,m+1},{\mathcal{O}})$ has a dense image if $p\geq n-1$.
Hence $H^{p}(D_{\varepsilon,1}\cup
D_{\varepsilon,2}\cup\cdots\cup D_{\varepsilon,m+1},{\mathcal{O}})$ vanishes for $p\geq n-1$ and
so is $H^{p}(D_{\varepsilon},{\mathcal{O}})$ if $p\geq \tilde{q}-2$, since $\tilde{q}-2+m=n-1$.\\
\hspace*{.1in}For the proof of the second assertion of the lemma,
we can write for all $0\leq j\leq k$ and $1\leq l\leq k$, $D_{j}\cap U_{l}=D'_{1}\cap\cdots\cap D'_{m+1}$,
where $D'_{i}=\{z\in U_{l}: \phi_{i,j}(z)<-\varepsilon\}$ are $q$-complete and $q$-Runge in $U_{l}$.
Therefore for all $t\leq m$ and $i_{1},\cdots,i_{t}\in\{1,\cdots,m+1\}$, $D'_{i_{1}}\cap
D'_{i_{2}}\cap\cdots\cap D'_{i_{t}}$ is $(\tilde{q}-1)$-Runge in $U_{l}$, which implies that
$H^{p}(U_{l},{\mathcal{O}})$ has a dense image in $H^{p}(D'_{i_{1}}\cap
D'_{i_{2}}\cap\cdots\cap D'_{i_{t}},{\mathcal{O}})$ for all $p\geq \tilde{q}-2$. This shows that
$H^{p}(D'_{i_{1}}\cap D'_{i_{2}}\cap\cdots\cap D'_{i_{t}},{\mathcal{O}})=0$ for all $p\geq \tilde{q}-2$. Moreover, by (~\cite{ref9}, Proposition $1$), one obtains
$$H^{p}(D_{j}\cap U_{l},{\mathcal{O}})\cong H^{p+m}(D'_{1}\cup\cdots\cup D'_{m+1},{\mathcal{O}}) \ \ for \ \ p\geq \tilde{q}-2.$$
Since $\tilde{q}-2+m=n-1$ and $D'_{1}\cup\cdots\cup D'_{m+1}$ is $n$-Runge in $U_{l}$, then
$$H^{p}(D_{j}\cap U_{l},{\mathcal{O}})\cong H^{p+m}(D'_{1}\cup\cdots\cup D'_{m+1},{\mathcal{O}})=0 \ \ for \ \ p\geq \tilde{q}-2.$$
We now consider the Mayer-Vietoris sequence for cohomology
\begin{center}
$\cdots\rightarrow H^{p-1}(D_{j}\cap U_{j+1},{\mathcal{O}})
\rightarrow H^{p}(D_{j+1},{\mathcal{O}})\rightarrow H^{p}(D_{j},
{\mathcal{O}})\oplus H^{p}(D_{j+1}\cap U_{j+1},
{\mathcal{O}})\rightarrow H^{p}(D_{j}\cap U_{j+1},
{\mathcal{O}})\rightarrow\cdots$
\end{center}
Since $H^{p}(D_{j}\cap
U_{j+1},{\mathcal{O}})=H^{p}(D_{j+1}\cap U_{j+1},{\mathcal{O}})=0$
for $p\geq \tilde{q}-2$, then $H^{p}(D_{j+1},{\mathcal{O}})\rightarrow
H^{p}(D_{j},{\mathcal{O}})$ is surjective for all $p\geq \tilde{q}-2$.\\
We are now going to show that
$H^{p}(D_{j+1},{\mathcal{O}})\rightarrow
H^{p}(D_{j},{\mathcal{O}})$ is injective for all $p\geq \tilde{q}-2$.
Let ${\mathcal{V}}=(V_{i})_{i\in \mathbb{N}}$ be an open
covering of $D_{\varepsilon_{o}}$ with a fundamental system of Stein neighborhoods of $D_{\varepsilon_{o}}$
such that if $V_{i_{o}}\cap\cdots \cap V_{i_{r}}\neq\emptyset$ and
$V_{i_{o}}\cup\cdots \cup V_{i_{r}}\subset D_{j+1}$, then
$V_{i_{o}}\cup\cdots \cup V_{i_{r}}\subset D_{j}$ or
$V_{i_{o}}\cup\cdots \cup V_{i_{r}}\subset U_{j+1}\cap D_{j+1}$.\\
\hspace*{.1in}We first show that $H^{n-1}(D_{k}, {\mathcal{O}})=0$. We shall prove it
assuming that it has already been proved for $j<k$. For this, we consider
the Mayer-Vietoris sequence for cohomology
$$\rightarrow
H^{\tilde{q}-3}(D_{j}, {\mathcal{O}})\oplus H^{\tilde{q}-3}(D_{j+1}\cap U_{j+1}, {\mathcal{O}})\stackrel{r^{*}}\rightarrow
H^{\tilde{q}-3}(D_{j}\cap U_{j+1}, {\mathcal{O}})\stackrel{j^{*}}\rightarrow H^{\tilde{q}-2}(D_{j+1}, {\mathcal{O}})
\stackrel{\rho^{*}}\rightarrow$$
Let $\xi$ be a cocycle in $Z^{\tilde{q}-2}({\mathcal{V}}|_{D_{j+1}}, {\mathcal{O}})$ and let
$\rho(\xi)$ be its restriction to a cocycle in $Z^{\tilde{q}-2}({\mathcal{V}}|_{D_{j}}, {\mathcal{O}})$.
Since $\rho(\xi)$ is a coboundary by induction and\\
$H^{\tilde{q}-2}(D_{j+1}\cap U_{j+1}, {\mathcal{O}})=0$, from the Mayer-Vietoris sequence,
it follows that there exist
$$\eta\in Z^{\tilde{q}-3}({\mathcal{V}}|_{D_{j}\cap U_{j+1}}, {\mathcal{O}}) \ \
and \ \ \mu\in C^{\tilde{q}-3}({\mathcal{V}}|_{D_{j+1}}, {\mathcal{O}})$$
such that $\xi=j(\eta)+\delta\mu$. There exists a sequence
$\{\eta_{n}\}\subset Z^{\tilde{q}-3}({\mathcal{V}}|_{D_{j+1}\cap U_{j+1}}, {\mathcal{O}})$
with $r(\eta_{n})-\eta\rightarrow 0$, when $n\rightarrow \infty$. This is possible
because\\
$Z^{\tilde{q}-3}({\mathcal{V}}|_{D_{j+1}\cap U_{j+1}}, {\mathcal{O}})\rightarrow Z^{\tilde{q}-3}({\mathcal{V}}|_{D_{j}\cap U_{j+1}}, {\mathcal{O}})$
has a dense range. (See the proof of lemma $2$). Now choose a sequence
$\{\gamma_{n}\}\subset C^{\tilde{q}-3}({\mathcal{V}}|{D_{j+1}}, {\mathcal{O}})$ such that
$j(r(\eta_{n}))=\delta\gamma_{n}$. Then
$$\xi-\delta\mu-\delta\gamma_{n}=j(\eta-r(\eta_{n}))$$
This proves that $\delta\mu+\delta\gamma_{n}$ converges to $\xi$ when $n\rightarrow \infty$.
Since, by lemma $1$, $dim_{\mathbb{C}}H^{\tilde{q}-2}(D_{j+1}, {\mathcal{O}})<\infty$, then the coboundary
space $B^{\tilde{q}-2}({\mathcal{V}}|_{D_{j+1}}, {\mathcal{O}})$ is closed in
$Z^{\tilde{q}-2}({\mathcal{V}}|_{D_{j+1}}, {\mathcal{O}})$. Therefore
$\xi\in B^{\tilde{q}-2}({\mathcal{V}}|_{D_{j+1}}, {\mathcal{O}})$ and $H^{\tilde{q}-2}(D_{j}, {\mathcal{O}})=0$
for all $0\leq j\leq k$.\\
\hspace*{.1in}Now since $D_{\varepsilon}\subset\subset D_{k}(\varepsilon)=D_{k}$, there exists $0<\varepsilon'<\varepsilon$ such that $\varepsilon'>\varepsilon_{0}$ and
$D_{\varepsilon'}\subset\subset D_{k}(\varepsilon)$.
Then we have $D_{\varepsilon}\subset D_{\varepsilon'}\subset D_{k}(\varepsilon)\subset D_{k}(\varepsilon').$
Since for every $p\geq \tilde{q}-2$,
$H^{p}(D_{k}(\varepsilon'), {\mathcal{O}})\rightarrow H^{p}(D_{\varepsilon'}, {\mathcal{O}})$ is surjective by lemma $1$, then $H^{p}(D_{k}(\varepsilon), {\mathcal{O}})\rightarrow H^{p}(D_{\varepsilon'}, {\mathcal{O}})$ is also surjective for all $p\geq \tilde{q}-2$, which shows that $H^{p}(D_{\varepsilon'}, {\mathcal{O}})=0$ for $p\geq \tilde{q}-2$ and $\varepsilon'\in A$.
\end{proof}
\begin{Propo}The set $D_{\varepsilon_{o}}$ is cohomologially $(\tilde{q}-3)$-complete with respect to the structure sheaf ${\mathcal{O}}_{D_{\varepsilon_{0}}}$, which means that $H^{p}(D_{\varepsilon_{0}},{\mathcal{O}}_{D_{\varepsilon_{0}}})=0$
for all $p\geq \tilde{q}-3$.
\end{Propo}
\begin{proof}
\hspace*{.1in}In order to prove the Proposition, we consider the set
$A$ of all real numbers $\varepsilon\geq \varepsilon_{o}$ such that
$H^{p}(D_{\varepsilon}, {\mathcal{O}})=0$ for all $p\geq
\tilde{q}-2$, where $D_{\varepsilon}=\{z\in\mathbb{C}^{n}:
\rho(z)<-\varepsilon\}$. Then, by lemma $3$, the set $A$ is not
empty and open in $[+\varepsilon_{0},+\infty[.$ Moreover, if
$\varepsilon=Inf(A)$, there exists a decreasing sequence of real
numbers $\varepsilon_{j}\in A$, $j\in\mathbb{N}$, such that
$\varepsilon_{j}\rightarrow \varepsilon$. Since
$H^{p}(D_{\varepsilon_{j}},{\mathcal{O}})=0$ for all $p\geq
\tilde{q}-2$ and, by lemma $2$,
$H^{p}(D_{\varepsilon_{j+1}},{\mathcal{O}})\rightarrow
H^{p}(D_{\varepsilon_{j}},{\mathcal{O}})$ has a dense image for
$p\geq \tilde{q}-3$, then $H^{p}(D_{\varepsilon},{\mathcal{O}})=0$
for all $p\geq \tilde{q}-2$ (see ~\cite{ref2}, p. $250$).
Hence
$\varepsilon\in A.$
\hspace*{.1in}Suppose now that $\varepsilon>\varepsilon_{0}$. then
there exists $\varepsilon'\in A$ such that
$\varepsilon_{0}<\varepsilon'<\varepsilon$, which is a
contradiction. Therefore $\varepsilon_{0}=\varepsilon\in A$.
But, since
$H^{\tilde{q}-3}(D_{\varepsilon_{0}},{\mathcal{O}})=0$ according to the proof of lemma $2$, then $D_{\varepsilon_{0}}$ is cohomologically $(\tilde{q}-3)$-complete.
\end{proof}
\newpage
\noindent
\begin{Thm}Let $(n,q)$ be a pair of integers such that $n\geq 3,$
$1<q<n$, and $q$ does not divide $n$.
We put $m=[\frac{n}{q}]$ and $r=n-mq$.
Then there exists a
cohomologically $(\tilde{q}-3)$-complete domain $D$ with respect to $O_{D}$
in $\mathbb{C}^{n}$ such that $H_{n+\tilde{q}-3}(D,\mathbb{Z})\neq
0$, if $r=1$ and $m>q\geq 2$, where
$\tilde{q}=n-[\frac{n}{q}]+1$ and $[\frac{n}{q}]$ is the
integral part of $\frac{n}{q}$.
\end{Thm}
\begin{proof}
We have proved in the Proposition that
$D_{\varepsilon_{o}}$ is cohomologically $(\tilde{q}-3)$-complete
with respect to the structure sheaf ${\mathcal{O}}_{D_{\varepsilon_{0}}}$.\\
\hspace*{.1in}It was shown by Diederich-Fornaess
~\cite{ref5} that if $\delta>0$ is small enough, then the
following topological sphere of real dimension $n+\tilde{q}-2$:
$$S_{\delta}=\{z\in\mathbb{C}^{n}:\,\,
x_{1}^{2}+\cdots+x_{m}^{2}+|z_{m+1}|^{2}+\cdots+|z_{n}|^{2}=\delta,$$
$$y_{j}=-\displaystyle\sum_{i=m+1}^{n}|z_{i}|^{2}+(m+1)(\displaystyle
\sum_{i=m+(j-1)(q-1)+1}^{m+j(q-1)}|z_{i}|^{2}) \mbox{\;for\;}
\hspace{1pc}$$ $$j=1,\cdots, m\},\hspace{14.7pc}$$ where
$z_{j}=x_{j}+y_{j}$ for $j=1,\cdots,n$, is not homologous to $0$
in $D_{\varepsilon_{o}}$. This follows from the fact that the set
$E=\{z\in \mathbb{C}^{n}:
x_{1}=\cdots=x_{m}=z_{m+1}=\cdots=z_{n}=0\}$ does not intersect
$D_{\varepsilon_{o}}$.\\
\hspace*{.1in}Let $0<\varepsilon<\varepsilon_{0}$ be such that $D_{\varepsilon}\subset\subset B$.
Then it follows, exactly
as for $D_{\varepsilon_{o}}$, that $D_{\varepsilon}$ is
cohomologically $(\tilde{q}-2)$-complete with respect to ${\mathcal{O}}_{D_{\varepsilon}}$.
We may take $\varepsilon$ such that $D_{\varepsilon}$ does not intersect $E$.
Then $S_{\delta}$ is not homologous to $0$ in $D_{\varepsilon}$.\\
\hspace*{.1in}We are now going to prove that $H_{n+\tilde{q}-3}(D_{\varepsilon},\mathbb{Z})$ does not vanish.
Indeed, we define for $1\leq j\leq m+1$ the sets
$E_{j}=(D_{1}\cap\cdots \cap D_{m-j+2})\cup (D_{m-j+3}\cup\cdots\cup D_{m+1}$).
Then by using the Mayer-Vietoris sequence for homology
\begin{center}
$0=H_{n+\tilde{q}+j-3}(D_{1}\cap\cdots\cap D_{m-j+2},\mathbb{Z})\oplus
H_{n+\tilde{q}+j-3}(D_{m-j+3}\cup\cdots\cup D_{m+1},\mathbb{Z})\rightarrow
H_{n+\tilde{q}+j-3}(E_{j+1},\mathbb{Z})\rightarrow H_{n+\tilde{q}+j-4}(E_{j},\mathbb{Z})
\rightarrow H_{n+\tilde{q}+j-4}(D_{1}\cap\cdots\cap D_{m-j+2},\mathbb{Z})\oplus H_{n+\tilde{q}+j-4}(D_{m-j+3}\cup\cdots\cup D_{m+1},\mathbb{Z})\rightarrow\cdots$
\end{center}
one can easily verify by induction that
$$H_{2n-2}(D_{1}\cup\cdots\cup D_{m+1},\mathbb{Z})\rightarrow H_{n+\tilde{q}-3}(D_{\varepsilon},\mathbb{Z})$$
is injective.\\
\hspace*{.1in}We first show that $H^{n+\tilde{q}-3}(D_{\varepsilon},\mathbb{C})\neq 0.$
For this, we consider the $(2n-m-2)$-real form defined as follows: \\
$\omega=\frac{1}{2}(\displaystyle\sum_{i=1}^{n}x_{i}^{2}+\displaystyle\sum_{i=m+2}^{n}y_{i}^{2})^{-2n+m+1}(\displaystyle\sum_{i=1}^{n}
(-1)^{i}x_{i}dx_{1}\wedge\cdots\wedge \widehat{dx_{i}}\wedge\cdots\wedge dx_{n}\wedge dy_{m+2}\wedge\cdots\wedge dy_{m+i}\wedge\cdots\wedge dy_{n}
+\displaystyle\sum_{i=2}^{n-m}(-1)^{n+i-1}
y_{m+i}dx_{1}\wedge\cdots\wedge dx_{n}\wedge dy_{m+2}\wedge\cdots\wedge\widehat{dy_{m+i}}\wedge\cdots\wedge dy_{n})$
\newpage
\noindent
Then $\omega$ is d-closed and therefore defines a cohomology class in $H^{2n-m-2}(D_{\varepsilon},\mathbb{C}).$
\noindent
\hspace*{.1in}Let $S'_{\delta}$ be the topological sphere of real dimension $2n-m-2$ defined by\\
$S'_{\delta}=\{z\in S_{\delta}: y_{m+1}=0\},$ where $z_{j}=x_{j}+iy_{j}$ for $j=1,\cdots, n$.\\
\hspace*{.1in}Since $\omega$ does not depend on $y_{1},\cdots, y_{m+1}$, then
$$\int_{S'_{\delta}}\omega\neq 0$$
This implies that $H^{2n-m-2}(D_{\varepsilon},\mathbb{C})\neq 0.$\\
\hspace*{.1in}We shall prove now that if $E=D_{1}\cup\cdots\cup D_{m+1}$, then the homology group $H_{2n-2}(E, \mathbb{Z})$ does not vanish.
We note first that the map
$$H^{2n-m-2}(D_{\varepsilon},\mathbb{C})\rightarrow H^{2n-2}(E,\mathbb{C})$$
is an isomorphism.
In fact, we consider over $E$ the resolution of the constant sheaf $\mathbb{C}$:
$$0\rightarrow \mathbb{C}\stackrel{i}\rightarrow {\mathcal{O}}\stackrel{d}\rightarrow\Omega^{1}
\stackrel{d}\rightarrow\cdots\stackrel{d}\rightarrow\Omega^{n}\rightarrow 0$$
where $\Omega^{p}$ is the sheaf of germs of holomorphic $p$-forms on $E$.
This resolution maybe break it up into short exact sequences
$$0\rightarrow Z^{i}\rightarrow \Omega^{i}\rightarrow Z^{i+1}\rightarrow 0, \ \ for \ \ i=0,\cdots, n$$
where $Z^{0}=\mathbb{C}$, $\Omega^{0}={\mathcal{O}}$, and $Z^{i}=Im\{\Omega^{i-1}\stackrel{d}\rightarrow \Omega^{i}\}$ for $1\leq i\leq n$.\\
\hspace*{.1in}On the other hand, since for any integer $t$ with $1\leq t\leq m-1$ and all $i_{1},\cdots, i_{t}\in\{1,\cdots, m+1\},$ $D_{i_{1}}\cap\cdots\cap D_{i_{t}}$ is $((m-1)(q-1)+1)$-Runge in $B$ and $n-m-2\geq (m-1)(q-1)$, then
$H^{p}(D_{i_{1}}\cap\cdots\cap D_{i_{t}},\mathbb{C})=0$ for $p\geq 2n-m-2$. Hence for any indexes
$j_{1},\cdots, j_{m}\in\{1,\cdots, m+1\}$ we have
$$H^{2n-m-2}(D_{j_{1}}\cap\cdots\cap D_{j_{m}},\mathbb{C})\cong H^{2n-3}(D_{j_{1}}\cup\cdots\cup D_{j_{m}},\mathbb{C})$$
But $H^{2n-3}(D_{j_{1}}\cup\cdots\cup D_{j_{m}},\mathbb{C})=0$. Indeed, if we put
$E_{m}=D_{j_{1}}\cup\cdots\cup D_{j_{m}},$ then obviously $H^{p}(E_{m},\Omega^{j})=0$ for all $p\geq n-1,$
and $j\geq 0$, because $E_{m}$ is $n$-Runge in $B$. Note also that since $n-m-2\geq (m-1)q-(m-1)$, then, by the proof of lemma $2$, $H^{n-2}(E_{m},\Omega^{n-1})\cong H^{n-m-1}(D_{j_{1}}\cap\cdots\cap D_{j_{m}},\Omega^{n-1})=0$ and
$H^{n-3}(E_{m},\Omega^{n})\cong H^{n-m-2}(D_{j_{1}}\cap\cdots\cap D_{j_{m}},\Omega^{n})=0$,
because the $\Omega^{j}$ are free and $n-m-2\geq q$ and $n-m-1\leq n-q-2$. (See the proof of lemma $2$).
Therefore from the long exact sequences of cohomology associated to the short exact sequences
$$0\rightarrow Z^{i}\rightarrow \Omega^{i}\rightarrow Z^{i+1}\rightarrow 0, \ \ for \ \ i=0,\cdots, n$$
we deduce a natural $\mathbb{C}$-isomorphism
$$H^{2n-3}(E_{m},\mathbb{C})\cong H^{n-2}(E_{m},Z^{n-1})$$
and, an exact sequence
$$\cdots\rightarrow H^{n-3}(E_{m},\Omega^{n})\rightarrow H^{n-2}(E_{m},Z^{n-1})
\rightarrow H^{n-2}(E_{m},\Omega^{n-1})\rightarrow\cdots$$
Since $H^{n-3}(E_{m},\Omega^{n})=H^{n-2}(E_{m},\Omega^{n-1})=0,$ then
$H^{2n-3}(E_{m},\mathbb{C})\cong H^{n-2}(E_{\varepsilon},Z^{n-1})=0$, which shows that
$H^{2n-m-2}(D_{\varepsilon},\mathbb{C})\rightarrow H^{2n-2}(E,\mathbb{C})$
is an isomorphism.\\
\hspace*{.1in}Now by using the universal coefficient theorem for homology
$$H_{2n-2}(E,\mathbb{C})\cong H_{2n-2}(E,\mathbb{Z})\otimes
\mathbb{C}\oplus Tor(H_{2n-3}(E,\mathbb{Z}),\mathbb{C})$$
and the fact that $Tor(H_{2n-3}(E,\mathbb{Z}),\mathbb{C})=0$,
we conclude that $H_{2n-2}(E,\mathbb{Z})\neq 0$.
Since we know already that the map $H_{2n-2}(E,\mathbb{Z})\rightarrow H_{n+\tilde{q}-3}(D_{\varepsilon},\mathbb{Z})$
is injective, then $H_{n+\tilde{q}-3}(D_{\varepsilon},\mathbb{Z})$ does not vanish.
This completes the proof of the theorem.

\end{proof}

\newpage
\noindent
Mohammed Mou\c{c}ouf\\
D\'epartement de math\'ematiques et Informatique,\\
Universit\'e Chouaib
Doukkali,\\
Facult\'e des Sciences,\\
B.P.20; 24000, El
Jadida,Morocco.\\
Email : {\bf moucouf@hotmail.com}\\\\
Youssef Alaoui\\
D\'epartement de math\'ematiques et Informatique,\\
Universit\'e Chouaib
Doukkali,\\
Facult\'e des Sciences,\\
B.P.20; 24000, El
Jadida,Morocco.\\
Email : {\bf y.alaoui@iav.ac.ma or comp5123ster@gmail.com}
\end{document}